\documentclass{article}
\usepackage{amssymb}
\usepackage{amsmath}
\usepackage[english]{babel}
\newcommand{\F}{\text{$\mathbb{F}_q$}}
\author{E.Yu. Lerner}
\title{
Tables of graphs of binary and ternary sequences differentiation}
\begin{document}
\maketitle

In~\cite{ar1}--\cite{ar3} V.I.~Arnold studied the following
dynamical system generated by the operator of the finite
differen\-ti\-a\-ti\-on. Let $x$ be a cyclic sequence of $n$
elements of the finite field $\mathbb{F}_q$ (the first element
immediately follows the $n$-th one). Let $M$ be the set of such
sequences. Let us define the operation $\Delta$ as the transition
from $x$ to the sequence of differences of the neighbouring
elements from $x$. The dynamic system $\Delta$ is set by the
oriented graph~$G_n$ with vertices marked as $x$, $x\in M$. From
each vertex $x$ there comes exactly one edge (leading to $\Delta
x$). The attractors of the dynamic system $\Delta$ represent the
finite cycles of the length $m$, let us mark them as $O_m$.

Let $p$ be the characteristics of the field $\F$. Let $m$ be the
maximum power $p$, by which $n$ can be divided. The tree of
$h=p^m$ height leads to every point of every cycle. The edges of
the tree are oriented towards the root (located on a cycle). To
every internal vertex of the tree different from the root there
lead $q$ number of edges. $q-1$ number of edges leads to the root.
The total amount of tree vertices is therefore $q^h$. If $q = p$
Arnold marks this tree as ~$T_{p^h}$. Graph~$G_n$ expands into
connected components. A separate connected component was marked as
$(O_m*T_{p^h})$. Graph $G_n$ is the sum of several connected
components.

Sequence $x$ is called complicated one if vertex $x$ in graph
$G_n$  lead to the cycle of the most length. Arnold has set up
several hypotheses about the quota of complicated sequences and
concrete representatives of the most complicated ones. These
hypotheses are formulated more accurate and proved for some
partial cases in~\cite{Garber},~\cite{instead},~\cite{lerner2}.

The aim of this work is to give the results of calculations
of~$G_n$ for $q=2$, $n=2,\ldots,160$ and $q=3$, $n=2,\ldots,100$.
Graphs $G_n$, $n=2,\ldots,25$ and it's maximal component for
$n=26,\ldots,50$ had been calculated in~\cite{Karpen}. The results
of Karpenkov's calculation differ from my results for $n=15$ and
$n=17$.

The algorithm of calculations is rather simple. It is based on the
properties of polynomials algebra $\F [[z]]$ analysed by the
module $z^n-1$. A detailed description of the algorithm is given
in~\cite{instead}. That algorithm can be used in more common case.
The algorithm suitable for operator~$\Delta$ situation described
in paper of A.~Garber~\cite{Garber}
(see~\cite[Notation~3]{instead}). But it is not realized.

The program written on Mathematica is represented the Appendix.
The calculations had taken 20 minutes of computer time (computer
Intel Core 2 Duo).

The results of calculations of graph $q=2$, $n\le 300$ è $q=3$,
$n\le 150$ and program in Mathematica format is available
at~http://kek.ksu.ru/kek2/myArnold.htm. This calculations had
taken some hours of computer time. The algorithm would have been
polynomial one if there existed a polynomial algorithm of integer
factorization for the numbers of $q^n-1$ type. If $n>100$, then
with $n$-doubling the calculations time increase on average more
than tenfold. The tables enable us to define more precisely the
Arnold hypotheses and to prove them.

\begin{center}
\begin{tabular}{|r||r|r|}
\multicolumn{3}{c}{Table~1. Graph of operator $\Delta:
\mathbb{F}_2^n\to\mathbb{F}_2^n$}\\
\hline n &number of component  & connected components of  $G_n$ ($q=2$)\\
\hline \hline
2 & 1 &  $({O_1}*{T_4})$ \\
3 & 2 &  $ (O_3*T_2)+(O_1*T_2)$ \\
4 & 1 &  $({O_1}*{T_{16}})$ \\
5 & 2 &  $ (O_{15}*T_2)+(O_1*T_2)$ \\
6 & 4 & $ 2(O_6*T_4)+(O_3*T_4)+(O_1*T_4)$ \\
\hline
7 & 10 & $ 9(O_7*T_2)+(O_1*T_2)$ \\
8 & 1 & $({O_1}*{T_{256}})$ \\
9 & 6 & $ 4(O_{63}*T_2)+(O_3*T_2)+(O_1*T_2)$ \\
10 & 10 & $ 8(O_{30}*T_4)+(O_{15}*T_4)+(O_1*T_4)$ \\
11 & 4 & $ 3(O_{341}*T_2)+(O_1*T_2)$ \\
\hline
12 & 24 & $ 20(O_{12}*T_{16})+2(O_6*T_{16})+(O_3*T_{16})+(O_1*T_{16})$ \\
13 & 6 & $ 5(O_{819}*T_2)+(O_1*T_2)$ \\
14 & 298 & $ 288(O_{14}*T_4)+9(O_7*T_4)+(O_1*T_4)$ \\
15 & 1096 & $ 1091(O_{15}*T_2)+3(O_5*T_2)+(O_3*T_2)+(O_1*T_2)$ \\
16 & 1 & $({O_1}*{T_{65536}})$ \\
\hline
17 & 260 & $ 256(O_{255}*T_2)+3(O_{85}*T_2)+(O_1*T_2)$ \\
18 & 526 & $ 518(O_{126}*T_4)+4(O_{63}*T_4)+2(O_6*T_4)+(O_3*T_4)+(O_1*T_4)$ \\
19 & 28 & $ 27(O_{9709}*T_2)+(O_1*T_2)$ \\
20 & 1098 & $ 1088(O_{60}*T_{16})+8(O_{30}*T_{16})+(O_{15}*T_{16})+(O_1*T_{16})$ \\
21 & 16660 & $ 16640(O_{63}*T_2)+9(O_{21}*T_2)+9(O_7*T_2)+(O_3*T_2)+(O_1*T_2)$ \\
\hline
22 & 1540 & $ 1536(O_{682}*T_4)+3(O_{341}*T_4)+(O_1*T_4)$ \\
23 & 2050 & $ 2049(O_{2047}*T_2)+(O_1*T_2)$ \\
24 & 2744 & $ 2720(O_{24}*T_{256})+20(O_{12}*T_{256})+2(O_6*T_{256})+\phantom{111}$\\
 & & $+(O_3*T_{256})+(O_1*T_{256})$ \\
\hline
25 & 658 & $ 656(O_{25575}*T_2)+(O_{15}*T_2)+(O_1*T_2)$ \\
26 & 10246 & $ 10240(O_{1638}*T_4)+5(O_{819}*T_4)+(O_1*T_4)$ \\
27 & 4870 & $ 4864(O_{13797}*T_2)+4(O_{63}*T_2)+(O_3*T_2)+(O_1*T_2)$ \\
28 & 599338 & $ 599040(O_{28}*T_{16})+288(O_{14}*T_{16})+9(O_7*T_{16})+\phantom{111}$\\
 & & $+(O_1*T_{16})$ \\
\hline
29 & 566 & $ 565(O_{475107}*T_2)+(O_1*T_2)$ \\
30 & 8948416 & $ 8947294(O_{30}*T_4)+1091(O_{15}*T_4)+24(O_{10}*T_4)+\phantom{111111}$\\
 & & $+2(O_6*T_4)+3(O_5*T_4)+(O_3*T_4)+(O_1*T_4)$ \\
31 & 34636834 & $ 34636833(O_{31}*T_2)+(O_1*T_2)$ \\
\hline
32 & 1 & $({O_1}*{T_{4294967296}})$ \\
33 & 4198408 & $ 4198403(O_{1023}*T_2)+3(O_{341}*T_2)+(O_3*T_2)+(O_1*T_2)$ \\
34 & 8421892 & $ 8421248(O_{510}*T_4)+256(O_{255}*T_4)+384(O_{170}*T_4)+\phantom{111}$\\
 & & $+3(O_{85}*T_4)+(O_1*T_4)$ \\
 \hline
35 & 4195604 & $ 4195264(O_{4095}*T_2)+320(O_{819}*T_2)+9(O_{105}*T_2)+\phantom{1111}$\\
 & & $+(O_{15}*T_2)+9(O_7*T_2)+(O_1*T_2)$ \\
36 & 17043806 & $ 17043260(O_{252}*T_{16})+518(O_{126}*T_{16})+4(O_{63}*T_{16})+\phantom{11111}$\\
 & & $+20(O_{12}*T_{16})+2(O_6*T_{16})+(O_3*T_{16})+(O_1*T_{16})$\\
\hline
37 & 21256 & $ 21255(O_{3233097}*T_2)+(O_1*T_2)$ \\
38 & 3538972 & $ 3538944(O_{19418}*T_4)+27(O_{9709}*T_4)+(O_1*T_4)$ \\
39 & 67158052 & $ 67108860(O_{4095}*T_2)+49161(O_{1365}*T_2)+\phantom{11111111111}$\\
 & & $+20(O_{819}*T_2)+9(O_{455}*T_2)+(O_3*T_2)+(O_1*T_2)$ \\
\hline \hline
\end{tabular}

\begin{tabular}{|r||r|r|}
\multicolumn{3}{c}{Table~1. Graph of operator $\Delta:
\mathbb{F}_2^n\to\mathbb{F}_2^n$}\\
\hline n &number of component  & connected components of  $G_n$ ($q=2$)\\
\hline  \hline
40 & 35791946 & $ 35790848(O_{120}*T_{256})+1088(O_{60}*T_{256})+\phantom{11111111111111}$\\
 & & $+8(O_{30}*T_{256})+(O_{15}*T_{256})+(O_1*T_{256})$ \\
41 & 26214476 & $ 26214400(O_{41943}*T_2)+75(O_{13981}*T_2)+(O_1*T_2)$ \\
42 & 8726292328 & $ 8726273920(O_{126}*T_4)+16640(O_{63}*T_4)+\phantom{1111111111111}$\\
 & & $+1458(O_{42}*T_4)+9(O_{21}*T_4)+288(O_{14}*T_4)+9(O_7*T_4)+$\\
 & & $+2(O_6*T_4)+(O_3*T_4)+(O_1*T_4)$ \\
 \hline
43 & 805355524 & $ 805355523(O_{5461}*T_2)+(O_1*T_2)$ \\
44 & 806094340 & $ 806092800(O_{1364}*T_{16})+1536(O_{682}*T_{16})+\phantom{11111111111111}$\\
 & & $+3(O_{341}*T_{16})+(O_1*T_{16})$ \\
45 & 4296053112 & $ 4295999488(O_{4095}*T_2)+49104(O_{1365}*T_2)+\phantom{11111111111111}$\\
 & & $+144(O_{455}*T_2)+3276(O_{315}*T_2)+4(O_{63}*T_2)+\phantom{111111}$\\
 & & $+1091(O_{15}*T_2)+3(O_5*T_2)+(O_3*T_2)+(O_1*T_2)$ \\
 \hline
46 & 4297066498 & $ 4297064448(O_{4094}*T_4)+2049(O_{2047}*T_4)+(O_1*T_4)$ \\
47 & 8388610 & $ 8388609(O_{8388607}*T_2)+(O_1*T_2)$ \\
48 & 89479864 & $ 89477120(O_{48}*T_{65536})+2720(O_{24}*T_{65536})+20(O_{12}*T_{65536})+$\\
 & & $+2(O_6*T_{65536})+(O_3*T_{65536})+(O_1*T_{65536})$\\
 \hline
49 & 134217802 & $ 134217792(O_{2097151}*T_2)+9(O_7*T_2)+(O_1*T_2)$ \\
50 & 5502932434 & $ 5502931768(O_{51150}*T_4)+656(O_{25575}*T_4)+\phantom{1111111111111}$\\
 & & $+8(O_{30}*T_4)+(O_{15}*T_4)+(O_1*T_4)$ \\
51 & 4415293883912 & $4415293686531(O_{255}*T_2)+197379(O_{85}*T_2)+\phantom{111111111}$\\
& &$+(O_3*T_2)+(O_1*T_2)$ \\
\hline
52 & 85920327686 & $ 85920317440(O_{3276}*T_{16})+10240(O_{1638}*T_{16})+\phantom{11111111}$\\
 & & $+5(O_{819}*T_{16})+(O_1*T_{16})$ \\
53 & 1266206 & $ 1266205(O_{3556769739}*T_2)+(O_1*T_2)$ \\
54 & 163209382798 & $ 163209377408(O_{27594}*T_4)+4864(O_{13797}*T_4)+\phantom{11111111111}$\\
 & & $+518(O_{126}*T_4)+4(O_{63}*T_4)+2(O_6*T_4)+\phantom{11111}$\\
& &$ +(O_3*T_4)+(O_1*T_4)$\\
 \hline
55 & 19327371272 & $ 16106142720(O_{1048575}*T_2)+3221228544(O_{349525}*T_2)+\phantom{11}$\\
 & & $+3(O_{5115}*T_2)+3(O_{341}*T_2)+(O_{15}*T_2)+(O_1*T_2)$ \\
56 & 5026339169578 & $ 5026338570240(O_{56}*T_{256})+599040(O_{28}*T_{256})+\phantom{11111}$\\
 & & $+288(O_{14}*T_{256})+9(O_7*T_{256})+(O_1*T_{256})$ \\
\hline 57 & 2473910599736 & $
2473910599707(O_{29127}*T_2)+27(O_{9709}*T_2)+\phantom{1111}$\\
& & $+(O_3*T_2)+(O_1*T_2)$ \\
58 & 75833016886 & $ 75833016320(O_{950214}*T_4)+565(O_{475107}*T_4)+(O_1*T_4)$ \\
59 & 9099508 & $9099507(O_{31675383749}*T_2)+(O_1*T_2)$\\
\hline
60 & 1200959905108816 & $ 1200959896157116(O_{60}*T_{16})+8947294(O_{30}*T_{16})+\phantom{1111111}$\\
 & & $+3264(O_{20}*T_{16})+1091(O_{15}*T_{16})+20(O_{12}*T_{16})+\phantom{1111}$\\
 & & $+24(O_{10}*T_{16})+2(O_6*T_{16})+3(O_5*T_{16})+\phantom{11}$\\
 & & $+(O_3*T_{16})+(O_1*T_{16})$\\
 \hline
61 & 17602326 & $ 17602325(O_{65498251203}*T_2)+(O_1*T_2)$ \\
62 & 18595508156138530 & $
18595508121501696(O_{62}*T_4)+34636833(O_{31}*T_4)+$\\
& &$+(O_1*T_4)$ \\
63 & 73201365371896620 & $ 73201365371846652(O_{63}*T_2)+49929(O_{21}*T_2)+\phantom{1111111}$\\
 & & $+28(O_9*T_2)+9(O_7*T_2)+(O_3*T_2)+(O_1*T_2)$ \\
\hline
64 & 1 & $({O_1}*{T_{18446744073709551616}})$ \\
65 & 4504699407499318 & $ 4504699407499263(O_{4095}*T_2)+48(O_{1365}*T_2)+\phantom{11}$\\
 & & $+5(O_{819}*T_2)+(O_{15}*T_2)+(O_1*T_2)$ \\
66 & 9016003948194832 & $ 9016003943994886(O_{2046}*T_4)+4198403(O_{1023}*T_4)+\phantom{1111}$\\
 & & $+1536(O_{682}*T_4)+3(O_{341}*T_4)+2(O_6*T_4)+(O_3*T_4)+$\\
 & & $+(O_1*T_4)$ \\
\hline \hline
\end{tabular}

\begin{tabular}{|r||r|r|}
\multicolumn{3}{c}{Table~1. Graph of operator $\Delta:
\mathbb{F}_2^n\to\mathbb{F}_2^n$}\\
\hline n &number of component  & connected components of  $G_n$ ($q=2$)\\
\hline \hline
67 & 128207980 & $ 128207979(O_{575525617597}*T_2)+(O_1*T_2)$ \\
68 & 18085043222151684 & $ 18085043201097728(O_{1020}*T_{16})+\phantom{1111111111111111111}$\\
 & & $+8421248(O_{510}*T_{16})+12632064(O_{340}*T_{16})+\phantom{11111}$\\
 & & $+256(O_{255}*T_{16})+384(O_{170}*T_{16})+3(O_{85}*T_{16})+$\\
 & & $+(O_1*T_{16})$ \\
69 & 70368769347588 & $ 70368756760576(O_{4194303}*T_2)+12582912(O_{1398101}*T_2)+\phantom{}$\\
 & & $+2049(O_{6141}*T_2)+2049(O_{2047}*T_2)+(O_3*T_2)+(O_1*T_2)$ \\
\hline
70 & 36037595295659812 & $ 36037595249505824(O_{8190}*T_4)+4195264(O_{4095}*T_4)+\phantom{11}$\\
 & & $+41953120(O_{1638}*T_4)+320(O_{819}*T_4)+4968(O_{210}*T_4)+$\\
 & & $+9(O_{105}*T_4)+8(O_{30}*T_4)+(O_{15}*T_4)+$\\
 & & $+288(O_{14}*T_4)+9(O_7*T_4)+(O_1*T_4)$ \\
71 & 34359738370 & $ 34359738369(O_{34359738367}*T_2)+(O_1*T_2)$ \\
\hline
72 & 36600682694456286 & $ 36600682677409760(O_{504}*T_{256})+\phantom{111111111111111111}$\\
 & & $+17043260(O_{252}*T_{256})+518(O_{126}*T_{256})+\phantom{1111111}$\\
 & & $+4(O_{63}*T_{256})+2720(O_{24}*T_{256})+20(O_{12}*T_{256})+$\\
 & & $+2(O_6*T_{256})+(O_3*T_{256})+(O_1*T_{256})$ \\
73 & 9241421688590306824 & $9241421688590303232(O_{511}*T_2)+3591(O_{73}*T_2)+$\\
& &$+(O_1*T_2)$ \\
\hline 74 & 730316239033096 & $730316239011840(O_{6466194}*T_4)+\phantom{1111111}$\\
& &$+21255(O_{3233097}*T_4)+(O_1*T_4)$ \\
75 & 18014415690024040 & $ 18014415689351152(O_{1048575}*T_2)+48(O_{349525}*T_2)+\phantom{1111}$\\
 & & $+671744(O_{25575}*T_2)+1091(O_{15}*T_2)+\phantom{1111}$\\
 & & $+3(O_5*T_2)+(O_3*T_2)+(O_1*T_2)$ \\
76 & 121597653799010332 & $ 121597653795471360(O_{38836}*T_{16})+\phantom{11111111111111}$\\
 & & $+3538944(O_{19418}*T_{16})+27(O_{9709}*T_{16})+(O_1*T_{16})$ \\
\hline
77 & 70368744243880 & $ 70368744243136(O_{1073741823}*T_2)+704(O_{97612893}*T_2)+\phantom{1}$\\
 & & $+27(O_{2387}*T_2)+3(O_{341}*T_2)+9(O_7*T_2)+(O_1*T_2)$ \\
78 & 9225625486372376656 & $ 9225623836634873850(O_{8190}*T_4)+67108860(O_{4095}*T_4)+$\\
 & & $+1649670162450(O_{2730}*T_4)+163870(O_{1638}*T_4)+\phantom{11}$\\
 & & $+49161(O_{1365}*T_4)+18432(O_{910}*T_4)+20(O_{819}*T_4)+$\\
 & & $+9(O_{455}*T_4)+2(O_6*T_4)+(O_3*T_4)+(O_1*T_4)$ \\
79 & 549755813890 & $ 549755813889(O_{549755813887}*T_2)+(O_1*T_2)$ \\
\hline
80 & 76861433658352714 & $ 76861433622560768(O_{240}*T_{65536})+\phantom{111111111111111}$\\
 & & $+35790848(O_{120}*T_{65536})+1088(O_{60}*T_{65536})+\phantom{111}$\\
 & & $+8(O_{30}*T_{65536})+(O_{15}*T_{65536})+(O_1*T_{65536})$ \\
81 & 111199991632646 & $ 111199991627776(O_{10871635887}*T_2)+4864(O_{13797}*T_2)+\phantom{1}$\\
 & & $+4(O_{63}*T_2)+(O_3*T_2)+(O_1*T_2)$ \\
\hline
82 & 14411532551533363276 & $ 14411532551467827200(O_{83886}*T_4)+\phantom{1111111111111111}$\\
 & & $+26214400(O_{41943}*T_4)+39321600(O_{27962}*T_4)+\phantom{11}$\\
 & & $+75(O_{13981}*T_4)+(O_1*T_4)$ \\
83 & 26494256092 &
$26494256091(O_{182518930210733}*T_2)+(O_1*T_2)$ \\
\hline
84 & 4797324681014830869808 & $ 4797324681006053048320(O_{252}*T_{16})+\phantom{111111111111111}$\\
 & & $+8726273920(O_{126}*T_{16})+50930100(O_{84}*T_{16})+\phantom{11111}$\\
 & & $+16640(O_{63}*T_{16})+1458(O_{42}*T_{16})+\phantom{1111111111111}$\\
 & & $+599040(O_{28}*T_{16})+9(O_{21}*T_{16})+288(O_{14}*T_{16})+\phantom{11}$\\
 & & $+20(O_{12}*T_{16})+9(O_7*T_{16})+2(O_6*T_{16})+\phantom{1}$\\
 & & $+(O_3*T_{16})+(O_1*T_{16})$ \\
\hline \hline
\end{tabular}

\begin{tabular}{|r||r|r|}
\multicolumn{3}{c}{Table~1. Graph of operator $\Delta:
\mathbb{F}_2^n\to\mathbb{F}_2^n$}\\
\hline n &number of component  &  connected components of  $G_n$ ($q=2$)\\
\hline  \hline
85 & 75854169073859085866528 & $ 75854169073859085340175(O_{255}*T_2)+197376(O_{85}*T_2)+$\\
 & & $+328960(O_{51}*T_2)+15(O_{17}*T_2)+(O_{15}*T_2)+(O_1*T_2)$ \\
86 & 1770995524065048969220 & $ 1770995524064243613696(O_{10922}*T_4)+\phantom{1111111}$\\
 & & $+805355523(O_{5461}*T_4)+(O_1*T_4)$ \\
\hline
87 & 288230377225455832 & $ 288230377225453567(O_{268435455}*T_2)+3(O_{89478485}*T_2)+\phantom{}$\\
 & & $+2260(O_{475107}*T_2)+(O_3*T_2)+(O_1*T_2)$ \\
88 & 443154625958844827140 & $ 443154625958038732800(O_{2728}*T_{256})+\phantom{1111111111111111}$\\
 & & $+806092800(O_{1364}*T_{256})+1536(O_{682}*T_{256})+\phantom{11111111}$\\
 & & $+3(O_{341}*T_{256})+(O_1*T_{256})$ \\
89 & 151189550474521284184066 & $ 151189550474521284184065(O_{2047}*T_2)+(O_1*T_2)$ \\
\hline 90 & 37788157488447270999544 & $ 37788157486791547953152(O_{8190}*T_4)+\phantom{11111111111111}$\\
 & & $+4295999488(O_{4095}*T_4)+1649668497048(O_{2730}*T_4)+\phantom{}$\\
 & & $+49104(O_{1365}*T_4)+4719672(O_{910}*T_4)+\phantom{111111111}$\\
 & & $+1744828722(O_{630}*T_4)+144(O_{455}*T_4)+$\\
 & & $+3276(O_{315}*T_4)+518(O_{126}*T_4)+4(O_{63}*T_4)+$\\
 & & $+8947294(O_{30}*T_4)+1091(O_{15}*T_4)+24(O_{10}*T_4)+$\\
 & & $+2(O_6*T_4)+3(O_5*T_4)+(O_3*T_4)+(O_1*T_4)$ \\
91 & 302305259898749761881162 & $ 302305259898749761880000(O_{4095}*T_2)+\phantom{111}$\\
 & & $+320(O_{819}*T_2)+832(O_{315}*T_2)+9(O_7*T_2)+(O_1*T_2)$ \\
\hline
92 & 37797387618632469579778 & $ 37797387618628172513280(O_{8188}*T_{16})+\phantom{1111111111}$\\
 & & $+4297064448(O_{4094}*T_{16})+2049(O_{2047}*T_{16})+\phantom{1111}$\\
 & & $+(O_1*T_{16})$ \\
93 & 4840430261137607980288068 & $ 4840430261134306154905600(O_{1023}*T_2)+\phantom{11111111}$\\
 & & $+3301756108800(O_{341}*T_2)+34636833(O_{93}*T_2)+\phantom{}$\\
 & & $+34636833(O_{31}*T_2)+(O_3*T_2)+(O_1*T_2)$ \\
\hline
94 & 295147940363733303298 & $ 295147940363724914688(O_{16777214}*T_4)+\phantom{111111111111}$\\
 & & $+8388609(O_{8388607}*T_4)+(O_1*T_4)$ \\
95 & 864691128470863928 & $ 864691128466931712(O_{22906492245}*T_2)+\phantom{111111111111}$\\
 & & $+3932160(O_{4581298449}*T_2)+27(O_{145635}*T_2)+\phantom{11111}$\\
 & & $+27(O_{9709}*T_2)+(O_{15}*T_2)+(O_1*T_2)$ \\
96 & 192153584145881784 & $ 192153584056401920(O_{96}*T_{4294967296})+\phantom{1111111111111}$\\
 & & $+89477120(O_{48}*T_{4294967296})+2720(O_{24}*T_{4294967296})+\phantom{}$\\
 & & $+20(O_{12}*T_{4294967296})+2(O_6*T_{4294967296})+\phantom{}$\\
 & & $+(O_3*T_{4294967296})+(O_1*T_{4294967296})$ \\
\hline
97 & 48684193446852291300 & $ 48684193446851772416(O_{1627389855}*T_2)+\phantom{1111111}$\\
 & & $+518883(O_{542463285}*T_2)+(O_1*T_2)$ \\
98 & 18889474938682197674314 & $ 18889474938682063456224(O_{4194302}*T_4)+\phantom{1111111}$\\
 & & $+134217792(O_{2097151}*T_4)+288(O_{14}*T_4)+\phantom{111111}$\\
 & & $+9(O_7*T_4)+(O_1*T_4)$ \\
\hline
99 & 97693956705350389350424 & $ 97693956705350372556800(O_{3243933}*T_2)+\phantom{11111111111}$\\
 & & $+12595212(O_{21483}*T_2)+4198403(O_{1023}*T_2)+\phantom{1111111}$\\
 & & $+3(O_{341}*T_2)+4(O_{63}*T_2)+(O_3*T_2)+(O_1*T_2)$ \\
100 & 774468841781667830583250 & $ 774468841781662327649728(O_{102300}*T_{16})+\phantom{111111111}$\\
 & & $+5502931768(O_{51150}*T_{16})+656(O_{25575}*T_{16})+\phantom{11111}$\\
 & & $+1088(O_{60}*T_{16})+8(O_{30}*T_{16})+(O_{15}*T_{16})+\phantom{1111}$\\
 & & $+(O_1*T_{16})$
\\ \hline \hline
\end{tabular}

\begin{center}
\begin{tabular}{|r||r|r|}
\multicolumn{3}{c}{Table~1. Graph of operator $\Delta:
\mathbb{F}_2^n\to\mathbb{F}_2^n$}\\
\hline n &number of component  &connected components of  $G_n$ ($q=2$)\\
\hline \hline
101 & 33442571490376 & $({O_1}*{T_2} + 33442571490375*{O_{37905296863701641}}*{T_2})$ \\
102 & 24855894122122 & $ 2485589412212211753180071430(O_{510}*T_4)+\phantom{1111111111111111111}$\\
 &17824209013776 & $+4415293686531(O_{255}*T_4)+1655735058432(O_{170}*T_4)+\phantom{111}$\\
 & & $+197379(O_{85}*T_4)+2(O_6*T_4)+(O_3*T_4)+(O_1*T_4)$ \\
\hline
103 & 2251799813685250 & $({O_1}*{T_2} + 2251799813685249*{O_{2251799813685247}}*{T_2})$ \\
104 & 1209221039595003 & $ 12092210395949947515043840(O_{6552}*T_{256})+\phantom{111111111111111111}$\\
 &3435371526 & $+85920317440(O_{3276}*T_{256})+10240(O_{1638}*T_{256})+\phantom{111111}$\\
 & & $+5(O_{819}*T_{256})+(O_1*T_{256})$ \\
\hline
105 & 4952969378181116 & $ 4952969378181116094346952704(O_{4095}*T_2)+\phantom{111111111111111111}$\\
 &115840687824 & $+21480079360(O_{819}*T_2)+13628160(O_{315}*T_2)+\phantom{111111111111}$\\
 & & $+9819(O_{105}*T_2)+16640(O_{63}*T_2)+27(O_{35}*T_2)+\phantom{111111}$\\
 & & $+9(O_{21}*T_2)+1091(O_{15}*T_2)+9(O_7*T_2)+\phantom{}$\\
 & & $+3(O_5*T_2)+(O_3*T_2)+(O_1*T_2)$ \\
\hline
106 & 2851240183087330 & $ 2851240183087329443840(O_{7113539478}*T_4)+\phantom{11111111111111111111}$\\
 &710046 & $+1266205(O_{3556769739}*T_4)+(O_1*T_4)$ \\
107 & 84179432287300 & $({O_1}*{T_2} + 84179432287299*{O_{963770320257286037}}*{T_2})$ \\
\hline
108 & 36751485112074499 & $ 367514851120744834390507520(O_{55188}*T_{16})+\phantom{1111111111111111}$\\
 &7616933598 & $+163209377408(O_{27594}*T_{16})+4864(O_{13797}*T_{16})+\phantom{11111111111}$\\
 & & $+17043260(O_{252}*T_{16})+518(O_{126}*T_{16})+4(O_{63}*T_{16})+\phantom{111}$\\
 & & $+20(O_{12}*T_{16})+2(O_6*T_{16})+(O_3*T_{16})+(O_1*T_{16})$ \\
\hline
109 & 11357291391466767 & $ 11357291391466767080488960(O_{28573587}*T_2)+\phantom{111111111111111111}$\\
 &080510606 & $+21645(O_{3174843}*T_2)+(O_1*T_2)$ \\
110 & 15595157945730310 & $ 154138188998497236961190400(O_{2097150}*T_4)+\phantom{111111111111111111}$\\
 &4619319336 & $+16106142720(O_{1048575}*T_4)+\phantom{11111111111111111111111111}$\\
 & & $+1813390458805848330729984(O_{699050}*T_4)+\phantom{111111111111111111}$\\
 & & $+3221228544(O_{349525}*T_4)+26136(O_{10230}*T_4)+\phantom{11111}$\\
 & & $+3(O_{5115}*T_4)+1536(O_{682}*T_4)+3(O_{341}*T_4)+\phantom{}$\\
 & & $+8(O_{30}*T_4)+(O_{15}*T_4)+(O_1*T_4)$ \\
\hline
111 & 56668397796909643 & $ 56668397794435742564340(O_{22906492245}*T_2)+\phantom{111111111111111111}$\\
 &811912 & $+2473901162523(O_{7635497415}*T_2)+27(O_{2545165805}*T_2)+\phantom{11}$\\
 & & $+85020(O_{3233097}*T_2)+(O_3*T_2)+(O_1*T_2)$ \\
\hline
112 & 70739430816307695 & $ 707394308163071929630064640(O_{112}*T_{65536})+\phantom{1111111111111111}$\\
 &5969234218 & $+5026338570240(O_{56}*T_{65536})+599040(O_{28}*T_{65536})+\phantom{11}$\\
 & & $+288(O_{14}*T_{65536})+9(O_7*T_{65536})+(O_1*T_{65536})$ \\
\hline
113 & 28047079119542908 & $ 2804707911954290859741020160(O_{1851279}*T_2)+\phantom{11111111111111111}$\\
 &59741020596 & $+435(O_{617093}*T_2)+(O_1*T_2)$ \\
114 & 89132022840231189 & $ 89132022840231186810739753014(O_{58254}*T_4)+\phantom{111111111111111111}$\\
 &284653891696 & $+2473910599707(O_{29127}*T_4)+3538944(O_{19418}*T_4)+\phantom{111}$\\
 & & $+27(O_{9709}*T_4)+2(O_6*T_4)+(O_3*T_4)+(O_1*T_4)$ \\
\hline 115 & 1180591620717478 & $1180591620717478412288({O_{17592186044415}}*{T_2})+2049({O_{30705}}*{T_2})+$ \\
 &416388 & $ + 2049*({O_{2047}}*{T_2})+({O_{15}}*{T_2})+({O_1}*{T_2})$\\
116 & 2732172362507197 & $ 2732172362507197092693278720(O_{1900428}*T_{16})+\phantom{11111111111111}$\\
 &168526295606 & $+75833016320(O_{950214}*T_{16})+565(O_{475107}*T_{16})+(O_1*T_{16})$\\
\hline
117 & 20287362573029851 &
     $ 20287362573029851539766654255104(O_{4095}*T_2)+\phantom{11111111111111}$\\
 &540612969432372 & $+824835096573(O_{1365}*T_2)+21480080620(O_{819}*T_2)+\phantom{1111111111}$\\
 & & $+9(O_{455}*T_2)+60(O_{273}*T_2)+4(O_{63}*T_2)+(O_3*T_2)+(O_1*T_2)$\\
\hline \hline
\end{tabular}
\end{center}

\begin{center}
\begin{tabular}{|r||r|r|}
\multicolumn{3}{c}{Table~1. Graph of operator $\Delta:
\mathbb{F}_2^n\to\mathbb{F}_2^n$}\\
\hline n &number of component  &connected components of  $G_n$ ($q=2$)\\
\hline \hline
118 & 13113771627025670 & $ 1311377162702567038255104(O_{63350767498}*T_4)+\phantom{11111111}$\\
 &47354612 & $+9099507(O_{31675383749}*T_4)+(O_1*T_4)$ \\
119 & 19807041809157774 & $ 19807041809157775483471855552(O_{16777215}*T_2)+\phantom{111111111111111}$\\
 & 548884102634& $+5368709440(O_{3355443}*T_2)+458752(O_{2396745}*T_2)+\phantom{11111}$\\
 & & $+2304(O_{1785}*T_2)+27(O_{595}*T_2)+256(O_{255}*T_2)+\phantom{}$\\
 & & $+3(O_{85}*T_2)+9(O_7*T_2)+(O_1*T_2)$ \\
\hline
120 & 43269140487790230 & $ 43269140487790229637274149969376(O_{120}*T_{256})+\phantom{1111111111111}$\\
 &838234162453456 & $+1200959896157116(O_{60}*T_{256})+107372544(O_{40}*T_{256})+\phantom{111111}$\\
 & & $+8947294(O_{30}*T_{256})+2720(O_{24}*T_{256})+\phantom{1111111111111111}$\\
 & & $+3264(O_{20}*T_{256})+1091(O_{15}*T_{256})+20(O_{12}*T_{256})+\phantom{}$\\
 & & $+24(O_{10}*T_{256})+2(O_6*T_{256})+3(O_5*T_{256})+\phantom{}$\\
 & & $+(O_3*T_{256})+(O_1*T_{256})$ \\
\hline
121 & 304904860722471940 & $304904860722471936(O_{4359484439294640007}*T_2)+\phantom{111111111111}$\\
 & & $+3(O_{341}*T_2)+(O_1*T_2)$ \\
122 & 101470495117893588 & $ 10147049511789358848409600(O_{130996502406}*T_4)+\phantom{11111111111}$\\
 &66011926 & $+17602325(O_{65498251203}*T_4)+(O_1*T_4)$ \\
\hline
123 & 507060723663925451 & $ 5070607236616196063404000018432(O_{1048575}*T_2)+\phantom{11111111111111}$\\
 & 8785019740312& $+3298534883328(O_{349525}*T_2)+\phantom{111111111111111111111111111111}$\\
 & & $+23058452082364252160(O_{209715}*T_2)+15728640(O_{69905}*T_2)+\phantom{}$\\
 & & $+104857675(O_{41943}*T_2)+75(O_{13981}*T_2)+(O_3*T_2)+(O_1*T_2)$ \\
\hline
124 & 107195806111686764 & $ 10719580611168676385087786739302400(O_{124}*T_{16})+\phantom{111111111}$\\
 &03683294895440930 & $+18595508121501696(O_{62}*T_{16})+34636833(O_{31}*T_{16})+\phantom{}$\\
 & & $+(O_1*T_{16})$ \\
125 & 151115727451828 & $ 151115727451828781056(O_{140737488355327875}*T_2)+\phantom{111111111}$\\
 & 781714& $+656(O_{25575}*T_2)+(O_{15}*T_2)+(O_1*T_2)$ \\
\hline
126 & 1687908566076083648 & $ 168790856607608364776581157512281594(O_{126}*T_4)+\phantom{1111111111111}$\\
 &49782549063005348 & $+73201365371846652(O_{63}*T_4)+26178823218(O_{42}*T_4)+\phantom{111111111}$\\
 & & $+49929(O_{21}*T_4)+3626(O_{18}*T_4)+288(O_{14}*T_4)+28(O_9*T_4)+\phantom{}$\\
 & & $+9(O_7*T_4)+2(O_6*T_4)+(O_3*T_4)+(O_1*T_4)$ \\
\hline 127 & 6698471789782253217 & \\
    & 78296471322378370 & $ 669847178978225321778296471322378369(O_{127}*T_2)+(O_1*T_2)$ \\
128 & 1 & $ (T_{340282366920938463463374607431768211456}*O_1)$ \\
129 & 2077045516211551385 & $ 20770455162115513853590587989606403(O_{16383}*T_2)+\phantom{11111}$\\
 &3590588794961928 & $+805355523(O_{5461}*T_2)+(O_3*T_2)+(O_1*T_2)$ \\
\hline
130 & 415485185495651359 & $ 41548518549565135951761643877136384(O_{8190}*T_4)+\phantom{1111111111}$\\
 &56266343286219174 & $+4504699407499263(O_{4095}*T_4)+1573224(O_{2730}*T_4)+\phantom{111111}$\\
 & & $+10240(O_{1638}*T_4)+48(O_{1365}*T_4)+5(O_{819}*T_4)+\phantom{111}$\\
 & & $+8(O_{30}*T_4)+(O_{15}*T_4)+(O_1*T_4)$ \\
\hline
131 & 281629680514649644 & $281629680514649643(O_{4833046947311902523261}*{T_2})+({O_1}*{T_2})$\\
132 & 831579586805812471 & $ 83157958680581247176179829486448700(O_{4092}*T_{16})+\phantom{1111111111}$\\
 &85195834240736352 & $+9016003943994886(O_{2046}*T_{16})+806092800(O_{1364}*T_{16})+\phantom{1111}$\\
 & & $+4198403(O_{1023}*T_{16})+1536(O_{682}*T_{16})+3(O_{341}*T_{16})+\phantom{1}$\\
 & & $+20(O_{12}*T_{16})+2(O_6*T_{16})+(O_3*T_{16})+(O_1*T_{16})$ \\
 \hline
133 & 207692666626040573 & $ 20769266662604057386289138824380415(O_{262143}*T_2)+\phantom{1111111111}$\\
 & 86289138826217168& $+1835015(O_{37449}*T_2)+1728(O_{9709}*T_2)+9(O_7*T_2)+(O_1*T_2)$ \\
134 & 47300395536410572 & $ 4730039553641057291367088128(O_{1151051235194}*T_4)+\phantom{11111111111}$\\
 &91495296108 & $+128207979(O_{575525617597}*T_4)+(O_1*T_4)$ \\
\hline
 \hline
\end{tabular}
\end{center}

\begin{center}
\begin{tabular}{|r||r|r|}
\multicolumn{3}{c}{Table~1. Graph of operator $\Delta:
\mathbb{F}_2^n\to\mathbb{F}_2^n$}\\
\hline n &number of component  &connected components of  $G_n$ ($q=2$)\\
\hline \hline 135 & 316912650061669076 & $ 316912650061669018800417144832(O_{68719476735}*T_2)+\phantom{111111111}$\\
 & 724008971640& $+52776356806656(O_{22906492245}*T_2)+\phantom{11111111111111111111111}$\\
 & & $+603979776(O_{7635497415}*T_2)+5142330998784(O_{896805}*T_2)+\phantom{11}$\\
 & & $+3983616(O_{68985}*T_2)+4864(O_{13797}*T_2)+\phantom{11111111111111111}$\\
 & & $+4295999488(O_{4095}*T_2)+49104(O_{1365}*T_2)+144(O_{455}*T_2)+$\\
 & & $+3276(O_{315}*T_2)+4(O_{63}*T_2)+1091(O_{15}*T_2)+3(O_5*T_2)+$\\
 & & $+(O_3*T_2)+(O_1*T_2)$ \\
\hline
136 & 1668050818239894429 & $ 166805081823989442856118195729727488(O_{2040}*T_{256})+\phantom{11111}$\\
 &01330803759841796 & $+18085043201097728(O_{1020}*T_{256})+\phantom{111111111111111111111}$\\
 & & $+27127564807962624(O_{680}*T_{256})+8421248(O_{510}*T_{256})+\phantom{}$\\
 & & $+12632064(O_{340}*T_{256})+256(O_{255}*T_{256})+\phantom{1111111111}$\\
 & & $+384(O_{170}*T_{256})+3(O_{85}*T_{256})+(O_1*T_{256})$ \\
\hline
137 & 370116963591829599 & $ 37011696359182959935519457280(O_{2353642078071}*T_2)+\phantom{1111111}$\\
 & 35895658796& $+376201515(O_{784547359357}*T_2)+(O_1*T_2)$ \\
138 & 103845961929503979 & $ 10384596192950287230240510760189952(O_{8388606}*T_4)+\phantom{11111111}$\\
 &10801735842013192 & $+70368756760576(O_{4194303}*T_4)+\phantom{111111111111111111111111111}$\\
 & & $+110680490830530084864(O_{2796202}*T_4)+\phantom{11111111111111111111}$\\
 & & $+12582912(O_{1398101}*T_4)+21485326338(O_{12282}*T_4)+\phantom{111111111}$\\
 & & $+2049(O_{6141}*T_4)+4297064448(O_{4094}*T_4)+2049(O_{2047}*T_4)+$\\
 & & $+2(O_6*T_4)+(O_3*T_4)+(O_1*T_4)$ \\
\hline
139 & 4246732448623781668 & $ 4246732448623781667(O_{82051117639860085587829}*T_2)+(O_1*T_2)$ \\
140 & 5318210374344337403 & $ 5318210374344337401807471618647917568(O_{16380}*T_{16})+\phantom{111111}$\\
 & 285013026940807588& $+36037595249505824(O_{8190}*T_{16})+4195264(O_{4095}*T_{16})+\phantom{111}$\\
 & & $+1441503810378787840(O_{3276}*T_{16})+41953120(O_{1638}*T_{16})+$\\
 & & $+320(O_{819}*T_{16})+2617842240(O_{420}*T_{16})+\phantom{1111111111111}$\\
 & & $+4968(O_{210}*T_{16})+9(O_{105}*T_{16})+1088(O_{60}*T_{16})+\phantom{111}$\\
 & & $+8(O_{30}*T_{16})+599040(O_{28}*T_{16})+(O_{15}*T_{16})+$\\
 & & $+288(O_{14}*T_{16})+9(O_7*T_{16})+(O_1*T_{16})$ \\
\hline
141 & 19807040628566365873 & $ 19807040628566365873362698240(O_{70368744177663}*T_2)+\phantom{11111111}$\\
 & 379475460& $+8388609(O_{25165821}*T_2)+8388609(O_{8388607}*T_2)+(O_3*T_2)+\phantom{}$\\
 & & $+(O_1*T_2)$ \\
142 & 20282409604241966234 & $ 20282409604241966234305956937728(O_{68719476734}*T_4)+\phantom{111}$\\
 & 340316676098& $+34359738369(O_{34359738367}*T_4)+(O_1*T_4)$ \\
\hline
143 & 48357032784585167030 & $ 4835703278458516703019008(O_{1152921504606846975}*T_2)+\phantom{1111}$\\
 & 19032& $+15(O_{279279}*T_2)+5(O_{819}*T_2)+3(O_{341}*T_2)+(O_1*T_2)$ \\
144 & 33758171321521672964 & $ 337581713215216729608063356501749760(O_{1008}*T_{65536})+\phantom{1111}$\\
 &4664039285683166 & $+36600682677409760(O_{504}*T_{65536})+\phantom{111111111111111111111}$\\
 & & $+17043260(O_{252}*T_{65536})+518(O_{126}*T_{65536})+\phantom{1111111111}$\\
 & & $+4(O_{63}*T_{65536})+89477120(O_{48}*T_{65536})+\phantom{1111111111111}$\\
 & & $+2720(O_{24}*T_{65536})+20(O_{12}*T_{65536})+2(O_6*T_{65536})+\phantom{}$\\
 & & $+(O_3*T_{65536})+(O_1*T_{65536})$ \\
\hline
145 & 93876727552027745924 & $ 9387672755202774592475260604702000686(O_{2375535}*T_2)+\phantom{1111}$\\
 &75260695701621176 & $+90999619923(O_{791845}*T_2)+565(O_{475107}*T_2)+(O_{15}*T_2)+\phantom{}$\\
 & & $+(O_1*T_2)$ \\
146 & 218206900181317251874 & $ 21820690018131725187408019477656969805824(O_{1022}*T_4)+\phantom{111}$\\
 & 17260899346030792200& $+9241421688590303232(O_{511}*T_4)+470679552(O_{146}*T_4)+\phantom{}$\\
 & & $+3591(O_{73}*T_4)+(O_1*T_4)$ \\
\hline
147 & 202824096036562821102 & $ 20282409603656282109965678673920(O_{4398046511103}*T_2)+\phantom{11111}$\\
 &10223507604 & $+244276381440(O_{18874359}*T_2)+134217792(O_{6291453}*T_2)+\phantom{11}$\\
 & & $+134217792(O_{2097151}*T_2)+16640(O_{63}*T_2)+9(O_{21}*T_2)+\phantom{}$\\
 & & $+9(O_7*T_2)+(O_3*T_2)+(O_1*T_2)$ \\
\hline
 \hline
\end{tabular}
\end{center}

\begin{center}
\begin{tabular}{|r||r|r|}
\multicolumn{3}{c}{Table~1. Graph of operator $\Delta:
\mathbb{F}_2^n\to\mathbb{F}_2^n$}\\
\hline n &number of component  &connected components of  $G_n$ ($q=2$)\\
\hline \hline
148 & 172441046452755849434 & $ 1724410464527558494342343881591357440(O_{12932388}*T_{16})+\phantom{11}$\\
 &3074197830390536 & $+730316239011840(O_{6466194}*T_{16})+21255(O_{3233097}*T_{16})+\phantom{}$\\
 & & $+(O_1*T_{16})$ \\
149 & 126774939137440139966 & $ 126774939137440139965(O_{2814530423790308547362667}*T_2)+(O_1*T_2)$\\
\hline
150 & 170141345719900803597 & $ 170141345719900803597526313209954303880(O_{2097150}*T_4)+\phantom{111}$\\
 &550097868581088208 & $+18014415689351152(O_{1048575}*T_4)+402653544(O_{699050}*T_4)+\phantom{}$\\
 & & $+48(O_{349525}*T_4)+5770242525159424(O_{51150}*T_4)+\phantom{1111111111}$\\
 & & $+671744(O_{25575}*T_4)+8947294(O_{30}*T_4)+1091(O_{15}*T_4)+\phantom{11}$\\
 & & $+24(O_{10}*T_4)+2(O_6*T_4)+3(O_5*T_4)+(O_3*T_4)+(O_1*T_4)$ \\
\hline
151 & 4355747223444196542430 & $ 43557472234441965424307564606143227527169(O_{32767}*T_2)+$\\
 &7564606143227527170 & $+(O_1*T_2)$ \\
152 & 2871143423438384893081 & $ 287114342343838489308000256284961013760(O_{77672}*T_{256})+\phantom{}$\\
 &21853938760024092 & $+121597653795471360(O_{38836}*T_{256})+\phantom{1111111111111111}$\\
 & & $+3538944(O_{19418}*T_{256})+27(O_{9709}*T_{256})+(O_1*T_{256})$ \\
\hline
153 & 3402823872033493138198 & $ 340282387203349257151469646223744434176(O_{16777215}*T_2)+\phantom{}$\\
 &70836020382667288 & $+56668401172135463092224(O_{5592405}*T_2)+\phantom{11111111111111}$\\
 & & $+13245881256972(O_{5355}*T_2)+4415293686531(O_{255}*T_2)+\phantom{}$\\
 & & $+197379(O_{85}*T_2)+4(O_{63}*T_2)+(O_3*T_2)+(O_1*T_2)$ \\
\hline
154 & 265845599404571182668 & $ 2658455994045711826684180298324736032(O_{2147483646}*T_4)+\phantom{}$\\
 &4274856325714408 & $+70368744243136(O_{1073741823}*T_4)+\phantom{111111111111111111111}$\\
 & & $+24189255833248(O_{195225786}*T_4)+704(O_{97612893}*T_4)+\phantom{11}$\\
 & & $+899424(O_{4774}*T_4)+27(O_{2387}*T_4)+1536(O_{682}*T_4)+\phantom{}$\\
 & & $+3(O_{341}*T_4)+288(O_{14}*T_4)+9(O_7*T_4)+(O_1*T_4)$ \\
\hline
155 & 217780922521473028604 & $ 21778092252147302860484537902010158022656(O_{1048575}*T_2)+\phantom{}$\\
 &84541279713169049668 & $+3377702941753344(O_{349525}*T_2)+34636833(O_{465}*T_2)+\phantom{}$\\
 & & $+34636833(O_{31}*T_2)+(O_{15}*T_2)+(O_1*T_2)$ \\
156 & 3485342350930329724696 & $ 348534235093029257742580260573643265556420(O_{16380}*T_{16})+\phantom{}$\\
 &23122057477492724232 & $+9225623836634873850(O_{8190}*T_{16})+\phantom{11111111111111111111}$\\
 & & $+3714727033635836197595726004(O_{5460}*T_{16})+\phantom{11111111111}$\\
 & & $+67108860(O_{4095}*T_{16})+21995602493740(O_{3276}*T_{16})+\phantom{1111}$\\
 & & $+1649670162450(O_{2730}*T_{16})+154656571392(O_{1820}*T_{16})+\phantom{}$\\
 & & $+163870(O_{1638}*T_{16})+49161(O_{1365}*T_{16})+\phantom{1111111111111}$\\
 & & $+18432(O_{910}*T_{16})+20(O_{819}*T_{16})+9(O_{455}*T_{16})+\phantom{}$\\
 & & $+20(O_{12}*T_{16})+2(O_6*T_{16})+(O_3*T_{16})+(O_1*T_{16})$ \\
\hline
157 & 86696145730328901894052& $ 8669614573032890189405275547333199285(O_{10536091491}*T_2)+\phantom{}$\\
 &75547333199286  & $+(O_1*T_2)$ \\
158 & 83076749736708357783940 & $ 83076749736708357783939769914359808(O_{1099511627774}*T_4)+\phantom{}$\\
 &319670173698 & $+549755813889(O_{549755813887}*T_4)+(O_1*T_4)$ \\
\hline
159 & 81129638414606699710187 & $ 81129638414606699710187514626047(O_{4503599627370495}*T_2)+\phantom{}$\\
 &519690872 & $+3(O_{1501199875790165}*T_2)+5064820(O_{3556769739}*T_2)+\phantom{}$\\
 & & $+(O_3*T_2)+(O_1*T_2)$ \\
160 & 7089215977519551322537 & $ 708921597751955132176933048662622208(O_{480}*T_{4294967296})+\phantom{}$\\
 &94482320974922 & $+76861433622560768(O_{240}*T_{4294967296})+\phantom{111111111111111}$\\
 & & $+35790848(O_{120}*T_{4294967296})+1088(O_{60}*T_{4294967296})+\phantom{}$\\
 & & $+8(O_{30}*T_{4294967296})+(O_{15}*T_{4294967296})+\phantom{}$\\
 & & $+(O_1*T_{4294967296})$ \\
 \hline \hline
\end{tabular}
\end{center}

\begin{tabular}{|r||r|r|}
\multicolumn{3}{c}{Table~2. Graph of operator $\Delta:
\mathbb{F}_3^n\to\mathbb{F}_3^n$}\\
\hline n &number of component  & connected components of  $G_n$ ($q=3$)\\
\hline \hline
2 & 3 & $3*({O_1}*{T_3})$ \\
3 & 1 & $ (O_1*T_{27})$ \\
4 & 6 & $3*({O_1}*{T_3}) + 3*({O_8}*{T_3})$ \\
5 & 2 & $ (O_{80}*T_3)+(O_1*T_3)$ \\
6 & 11 & $ 8(O_3*T_{27})+3(O_1*T_{27})$ \\
\hline
7 & 3 & $ 2(O_{364}*T_3)+(O_1*T_3)$ \\
8 & 279 & $3*({O_1}*{T_3}) + 6*({O_4}*{T_3}) + 270*({O_8}*{T_3})$ \\
9 & 1 & $ (O_1 T_{19683})$ \\
10 & 252 & $ 243(O_{80}*T_3)+6(O_{40}*T_3)+3(O_1*T_3)$ \\
11 & 246 & $ 243(O_{242}*T_3)+2(O_{121}*T_3)+(O_1*T_3)$ \\
\hline
12 & 832 & $ 818(O_{24}*T_{27})+3(O_8*T_{27})+8(O_3*T_{27})+3(O_1*T_{27})$ \\
13 & 20469 & $ 20412(O_{26}*T_3)+56(O_{13}*T_3)+(O_1*T_3)$ \\
14 & 4389 & $ 4374(O_{364}*T_3)+12(O_{182}*T_3)+3(O_1*T_3)$ \\
15 & 2216 & $ 2214(O_{240}*T_{27})+(O_{80}*T_{27})+(O_1*T_{27})$ \\
16 & 181800 & $3*({O_1}*{T_3}) + 6*({O_4}*{T_3}) + 270*({O_8}*{T_3}) +\phantom{11111111111111}$\\
 & & $+ 4374*({O_{40}}*{T_3}) + 177147*({O_{80}}*{T_3})$ \\
\hline
17 & 194 & $ 193(O_{223040}*T_3)+(O_1*T_3)$ \\
18 & 2195 & $ 2184(O_9*T_{19683})+8(O_3*T_{19683})+3(O_1*T_{19683})$ \\
19 & 519 & $ 518(O_{747916}*T_3)+(O_1*T_3)$ \\
20 & 14528301 & $ 14528241(O_{80}*T_3)+54(O_{40}*T_3)+3(O_8*T_3)+3(O_1*T_3)$ \\
21 & 354783 & $ 354780(O_{1092}*T_{27})+2(O_{364}*T_{27})+(O_1*T_{27})$ \\
\hline
22 & 43402482 & $ 43046721(O_{242}*T_3)+355758(O_{121}*T_3)+3(O_1*T_3)$ \\
23 & 177150 & $ 177147(O_{177146}*T_3)+2(O_{88573}*T_3)+(O_1*T_3)$ \\
24 & 435849063 & $ 435847140(O_{24}*T_{27})+1636(O_{12}*T_{27})+270(O_8*T_{27})+\phantom{11}$\\
 & & $+6(O_4*T_{27})+8(O_3*T_{27})+3(O_1*T_{27})$ \\
25 & 95663 & $ 95661(O_{2952400}*T_3)+(O_{80}*T_3)+(O_1*T_3)$ \\
\hline
26 & 32588084763 & $ 32587962120(O_{26}*T_3)+122640(O_{13}*T_3)+3(O_1*T_3)$ \\
27 & 1 & $ (O_1*T_{7625597484987})$ \\
28 & 10474724022 & $ 10474719630(O_{728}*T_3)+4374(O_{364}*T_3)+12(O_{182}*T_3)+\phantom{1111}$\\
 & & $+3(O_8*T_3)+3(O_1*T_3)$ \\
29 & 82466 & $ 82465(O_{277412144}*T_3)+(O_1*T_3)$ \\
\hline
30 & 31773382818 & $ 31773262986(O_{240}*T_{27})+119572(O_{120}*T_{27})+\phantom{111111111111}$\\
 & & $+243(O_{80}*T_{27})+6(O_{40}*T_{27})+8(O_3*T_{27})+3(O_1*T_{27})$ \\
31 & 231435 & $ 231434(O_{889632172}*T_3)+(O_1*T_3)$ \\
32 & 94172058441& $94143178827*({O_{6560}}*{T_3})+ 28697814*({O_{3280}}*{T_3}) +$ \\
 & & $+ 177147*({O_{80}}*{T_3})+ 4374*({O_{40}}*{T_3})+\phantom{11111}$\\
 & & $ 270*({O_8}*{T_3}) + 6*({O_4}*{T_3}) + 3*({O_1}*{T_3})$\\
\hline
33 & 283596620610 & $ 283596580836(O_{726}*T_{27})+39528(O_{363}*T_{27})+\phantom{111111111111}$\\
 & & $+243(O_{242}*T_{27})+2(O_{121}*T_{27})+(O_1*T_{27})$ \\
34 & 24924056094 & $ 24924051459(O_{223040}*T_3)+4632(O_{27880}*T_3)+3(O_1*T_3)$ \\
35 & 31381118670 & $ 31381118658(O_{531440}*T_3)+8(O_{7280}*T_3)+2(O_{364}*T_3)+\phantom{1111}$\\
 & & $+(O_{80}*T_3)+(O_1*T_3)$ \\
\hline
36 & 105911078650 & $ 105911075634(O_{72}*T_{19683})+818(O_{24}*T_{19683})+\phantom{1111111111111}$\\
 & & $+2184(O_9*T_{19683})+3(O_8*T_{19683})+8(O_3*T_{19683})+\phantom{1111}$\\
 & & $+3(O_1*T_{19683})$ \\
37 & 103053853533 & $ 103053850074(O_{1456468}*T_3)+3458(O_{112036}*T_3)+(O_1*T_3)$ \\
38 & 602051446125 & $ 602051439906(O_{747916}*T_3)+6216(O_{186979}*T_3)+3(O_1*T_3)$ \\
\hline \hline
\end{tabular}

\begin{tabular}{|r||r|r|}
\multicolumn{3}{c}{Table~2. Graph of operator $\Delta:
\mathbb{F}_3^n\to\mathbb{F}_3^n$}\\
\hline n &number of component  & connected components of  $G_n$ ($q=3$)\\
\hline \hline
39 & 1924290201095949 & $ 1924290191141640(O_{78}*T_{27})+9933840(O_{39}*T_{27})+\phantom{11111111}$\\
 & & $+20412(O_{26}*T_{27})+56(O_{13}*T_{27})+(O_1*T_{27})$ \\
40 & 50656939427274534 & $ 50656939398206748(O_{80}*T_3)+29056428(O_{40}*T_3)+\phantom{111111111111}$\\
 & & $+96(O_{20}*T_3)+10935(O_{16}*T_3)+270(O_8*T_3)+\phantom{1111111}$\\
 & & $+48(O_5*T_3)+6(O_4*T_3)+3(O_1*T_3)$ \\
\hline
41 & 1853302661441610 & $ 1853302661441604(O_{6560}*T_3)+5(O_{1312}*T_3)+(O_1*T_3)$ \\
42 & 3711131102090233 & $ 3711131082927684(O_{1092}*T_{27})+19158152(O_{546}*T_{27})+\phantom{11111}$\\
 & & $+4374(O_{364}*T_{27})+12(O_{182}*T_{27})+8(O_3*T_{27})+3(O_1*T_{27})$\\
43 & 121632015 & $ 121632014(O_{899590375372}*T_3)+(O_1*T_3)$ \\
\hline
44 & 5559154839586035 & $ 5559154709734350(O_{59048}*T_3)+86449200(O_{968}*T_3)+\phantom{11111111111}$\\
 & & $+43046721(O_{242}*T_3)+355758(O_{121}*T_3)+3(O_8*T_3)+3(O_1*T_3)$ \\
45 & 208464771247310 & $ 208464771245094(O_{720}*T_{19683})+2214(O_{240}*T_{19683})+\phantom{1111}$\\
 & & $+(O_{80}*T_{19683})+(O_1*T_{19683})$ \\
46 & 16677275843376846 & $ 16677275843376837(O_{177146}*T_3)+6(O_{88573}*T_3)+3(O_1*T_3)$ \\
\hline
47 & 94143178830 & $ 94143178827(O_{94143178826}*T_3)+2(O_{47071589413}*T_3)+(O_1*T_3)$ \\
48 & 12309659440396301526 & $ 12309613114542720714(O_{240}*T_{27})+\phantom{11111111111111111111111111}$\\
 & & $+46325417550228(O_{120}*T_{27})+177147(O_{80}*T_{27})+\phantom{1111111111}$\\
 & & $+1636(O_{12}*T_{27})+4374(O_{40}*T_{27})+435847140(O_{24}*T_{27})+\phantom{1}$\\
 & & $+270(O_8*T_{27})+6(O_4*T_{27})+8(O_3*T_{27})+3(O_1*T_{27})$\\
 & & $$ \\
\hline
49 & 77812219245 & $ 77812219242(O_{1025114613796}*T_3)+2(O_{364}*T_3)+(O_1*T_3)$ \\
50 & 81052475691459681 & $ 81052475690885463(O_{2952400}*T_3)+573966(O_{1476200}*T_3)+\phantom{1111}$\\
 & & $+243(O_{80}*T_3)+6(O_{40}*T_3)+3(O_1*T_3)$ \\
51 & 119210968252141016 & $ 119210968252140822(O_{669120}*T_{27})+193(O_{223040}*T_{27})+\phantom{11}$\\
 & & $+(O_1*T_{27})$ \\
\hline
52 & 2958371448703049336634 & $ 2958370209869331822030(O_{728}*T_3)+\phantom{11111111111111111111111}$\\
 & & $+1235346792567894(O_{364}*T_3)+3389148060480(O_{182}*T_3)+\phantom{11}$\\
 & & $+65176046880(O_{104}*T_3)+12754584(O_{91}*T_3)+\phantom{1111111111111}$\\
 & & $+32587962120(O_{26}*T_3)+122640(O_{13}*T_3)+3(O_8*T_3)+\phantom{}$\\
 & & $+3(O_1*T_3)$ \\
53 & 23979866306 & $ 23979866305(O_{269437777802768}*T_3)+(O_1*T_3)$ \\
\hline
54 & 282429537947 & $ 282429535752(O_{27}*T_{7625597484987})+\phantom{111111111111111111111111}$\\
 & & $+2184(O_9*T_{7625597484987})+8(O_3*T_{7625597484987})+\phantom{111111}$\\
 & & $+3(O_1*T_{7625597484987})$ \\
55 & 17294855100911367 & $ 16471290572295840(O_{3486784400}*T_3)+\phantom{1111111111111111111111}$\\
 & & $+823564528614792(O_{871696100}*T_3)+488(O_{9680}*T_3)+\phantom{111111}$\\
 & & $+243(O_{242}*T_3)+2(O_{121}*T_3)+(O_{80}*T_3)+(O_1*T_3)$ \\
\hline
56 & 239628037100450284331577 & $ 239628037100429331694620(O_{728}*T_3)+\phantom{11111111111111111111}$\\
 & & $+20946250614(O_{364}*T_3)+6386040(O_{182}*T_3)+24(O_{91}*T_3)+$\\
 & & $+270(O_8*T_3)+6(O_4*T_3)+3(O_1*T_3)$ \\
57 & 25916340428176453179 & $ 25916340428176452660(O_{2243748}*T_{27})+\phantom{1111111111111111111111}$\\
 & & $+518(O_{747916}*T_{27})+(O_1*T_{27})$ \\
\hline
58 & 5659604069395571388 & $ 5659604069395076595(O_{277412144}*T_3)+\phantom{111111111111111111111111}$\\
 & & $+494790(O_{138706072}*T_3)+3(O_1*T_3)$ \\
59 & 68630377364886 & $ 68630377364883(O_{68630377364882}*T_3)+\phantom{111111111111111111111111}$\\
 & & $+2(O_{34315188682441}*T_3)+(O_1*T_3)$ \\
\hline
60 & 6541845412842006768439897 & $ 6541845412842006666741342(O_{240}*T_{27})+\phantom{111111111111111111111}$\\
 & & $+87169428(O_{120}*T_{27})+14528241(O_{80}*T_{27})+54(O_{40}*T_{27})+$\\
 & & $+818(O_{24}*T_{27})+3(O_8*T_{27})+8(O_3*T_{27})+3(O_1*T_{27})$\\
\hline \hline
\end{tabular}

\begin{tabular}{|r||r|r|}
\multicolumn{3}{c}{Table~2. Graph of operator $\Delta:
\mathbb{F}_3^n\to\mathbb{F}_3^n$}\\
\hline n &number of component  & connected components of  $G_n$ ($q=3$)\\
\hline\hline
61 & 14358202911264125 &  \\
 &29274301 &$1435820291126412529274300*({O_{29524}}*{T_3})+({O_1}*{T_3})$\\
62 & 14295062479559704 & $ 142950624795578989998(O_{889632172}*T_3)+\phantom{11111111111111}$\\
 &1853 & $+18051852(O_{34216622}*T_3)+3(O_1*T_3)$ \\
\hline
63 & 17750224970402948
& $17750224970402948505180*({O_{3276}}*{T_{19683}})+\phantom{111111111111}$ \\
 &859963
& $+ 354780*({O_{1092}}*{T_{19683}}) +\phantom{111111111111}$\\
&
& $+ 2*({O_{364}}*{T_{19683}})+({O_1}*{T_{19683}})$\\
\hline
64 & 26588815594398467 & $ 26588814358957503287787(O_{43046720}*T_3)+\phantom{11111111111111111}$\\
 & 914122 & $+1235346792567894(O_{21523360}*T_3)+94143178827(O_{6560}*T_3)+$\\
 & & $+28697814(O_{3280}*T_3)+177147(O_{80}*T_3)+4374(O_{40}*T_3)+$\\
 & & $+270(O_8*T_3)+6(O_4*T_3)+3(O_1*T_3)$ \\
\hline
65 & 646109404691585649 & $ 6461094046914726761294802(O_{531440}*T_3)+\phantom{11111111111111111}$\\
 & 0130896 & $+1129718145924(O_{132860}*T_3)+10628820(O_{26572}*T_3)+\phantom{1111}$\\
 & & $+40880(O_{1040}*T_3)+(O_{80}*T_3)+20412(O_{26}*T_3)+\phantom{11111}$\\
 & & $+56(O_{13}*T_3)+(O_1*T_3)$ \\
\hline
66 & 157653078952793467 & $ 1576530569785667242428405684(O_{726}*T_{27})+\phantom{111111111111111111}$\\
 & 1872543206& $+219742267429400735032(O_{363}*T_{27})+43046721(O_{242}*T_{27})+\phantom{11}$\\
 & & $+355758(O_{121}*T_{27})+8(O_3*T_{27})+3(O_1*T_{27})$ \\
67 & 130186713272821981 &  \\
   & 0673103&$1301867132728219810673102*({O_{23737564}}*{T_3})+({O_1}*{T_3})$\\
\hline
68 & 215369416320196596 & $ 2153693963075556604049280(O_{43046720}*T_3)+\phantom{1111111111111111}$\\
 &5598801 & $+200126185045044696(O_{10761680}*T_3)+\phantom{1111111111111111}$\\
 & & $+224316463131(O_{223040}*T_3)+41688(O_{27880}*T_3)+\phantom{11}$\\
 & & $+3(O_8*T_3)+3(O_1*T_3)$ \\
\hline
69 & 581500652618604953 & $ 58150065261860474416260588(O_{531438}*T_{27})+\phantom{11111111111111111}$\\
 &37262242 & $+20920824504(O_{265719}*T_{27})+177147(O_{177146}*T_{27})+\phantom{11111111}$\\
 & & $+2(O_{88573}*T_{27})+(O_1*T_{27})$ \\
\hline
70 & 157004585340041586 & $ 1570045853400278602994617446(O_{531440}*T_3)+\phantom{1111111111111111}$\\
 & 4021980060& $+137261009820966(O_{265720}*T_3)+136080(O_{75920}*T_3)+\phantom{1111111}$\\
 & & $+15943200(O_{53144}*T_3)+3360(O_{37960}*T_3)+210(O_{7592}*T_3)+\phantom{}$\\
 & & $+1419120(O_{7280}*T_3)+35040(O_{3640}*T_3)+4374(O_{364}*T_3)+\phantom{}$\\
 & & $+12(O_{182}*T_3)+243(O_{80}*T_3)+6(O_{40}*T_3)+3(O_1*T_3)$ \\
\hline
71 & 50031545098999710 & $ 50031545098999707(O_{50031545098999706}*T_3)+\phantom{1111111111111}$\\
 & & $+2(O_{25015772549499853}*T_3)+(O_1*T_3)$ \\
\hline
72 & 158966843532060764 & $ 15896684353206076211804053860(O_{72}*T_{19683})+\phantom{111111111111111}$\\
 &24062056375 & $+211822151268(O_{36}*T_{19683})+435847140(O_{24}*T_{19683})+\phantom{111111}$\\
 & & $+1636(O_{12}*T_{19683})+2184(O_9*T_{19683})+270(O_8*T_{19683})+\phantom{}$\\
 & & $+6(O_4*T_{19683})+8(O_3*T_{19683})+3(O_1*T_{19683})$ \\
\hline
73 & 21195619020904687652 & $ 211956190209046876522656818030(O_{106288}*T_3)+(O_1*T_3)$ \\
   &  2656818031 & \\
74 & 464034902482014872 & $ 46403490248201486599936652622(O_{1456468}*T_3)+\phantom{111111111111}$\\
 & 18259766635& $+618323100444(O_{728234}*T_3)+3192(O_{364117}*T_3)+\phantom{11}$\\
 & & $+10374(O_{112036}*T_3)+3(O_1*T_3)$ \\
\hline\hline
\end{tabular}

\begin{tabular}{|r||r|r|}
\multicolumn{3}{c}{Table~2. Graph of operator $\Delta:
\mathbb{F}_3^n\to\mathbb{F}_3^n$}\\
\hline n &number of component  & connected components of  $G_n$ ($q=3$)\\
\hline\hline
75 & 2543512571121706003  & $2543512571121706003233512334*({O_{8857200}}*{T_{27}})+$ \\
 &233610211 & $+ 95661*({O_{2952400}}*{T_{27}}) +\phantom{11111111111}$\\
& &$+2214*({O_{240}}*{T_{27}})+ ({O_{80}}*{T_{27}})+({O_1}*{T_{27}})$\\
76 & 4066411119118709250 & $ 406641111911870924455019070810(O_{1495832}*T_3)+\phantom{1111}$\\
 &57070516938 & $+602051439906(O_{747916}*T_3)+6216(O_{186979}*T_3)+\phantom{1}$\\
 & & $+3(O_8*T_3)+3(O_1*T_3)$ \\
\hline
77 & 88629381196525444 & $ 8862938119652544083601(O_{205891132094648}*T_3)+\phantom{1111111111}$\\
 & 98168& $+413343(O_{29413018870664}*T_3)+976(O_{44044}*T_3)+\phantom{11111}$\\
 & & $+2(O_{364}*T_3)+243(O_{242}*T_3)+2(O_{121}*T_3)+(O_1*T_3)$ \\
\hline
78 & 7798292150171252732 & $ 7798292150171252629065743491190640(O_{78}*T_{27})+\phantom{111111111111}$\\
 & 977446669464691& $+103911670590189280(O_{39}*T_{27})+32587962120(O_{26}*T_{27})+\phantom{}$\\
 & & $+122640(O_{13}*T_{27})+8(O_3*T_{27})+3(O_1*T_{27})$ \\
\hline
79 & 25649083246955547 & $ 25649083246955546(O_{640303714176998250028}*T_3)+(O_1*T_3)$ \\
80 & 6158701225597746805
& $ 615870122559774680433023103113680536(O_{80}*T_3)+\phantom{11}$\\
 &34336981969109381 & $+101313878796417816(O_{40}*T_3)+58112088(O_{20}*T_3)+\phantom{}$\\
 & & $+896670(O_{16}*T_3)+1944(O_{10}*T_3)+270(O_8*T_3)+\phantom{}$\\
 & & $+48(O_5*T_3)+6(O_4*T_3)+3(O_1*T_3)$ \\
\hline
81 & 1 & $ (O_1*T_{443426488243037769948249630619149892803})$ \\
82 & 6759550125656063566 & $ 67595501256560635657283023512496812(O_{6560}*T_3)+\phantom{1111111111111}$\\
 & 8402839739466242 & $+11119815710329932(O_{3280}*T_3)+516639384(O_{1640}*T_3)+\phantom{}$\\
 & & $+15(O_{1312}*T_3)+96(O_{205}*T_3)+3(O_1*T_3)$ \\
83 & 36472996377170786406 & $ 36472996377170786403(O_{36472996377170786402}*T_3)+\phantom{111111111111}$\\
 & & $+2(O_{18236498188585393201}*T_3)+(O_1*T_3)$ \\
\hline
84 & 2030341063383872573
& $203034106338387257300456532894017860*({O_{2184}}*{T_{27}})+\phantom{11111}$ \\
 & 04167674470828544
& $+ 3711131082927684*({O_{1092}}*{T_{27}}) +\phantom{1111111111111111111}$\\
 &
& $+ 10474719630*({O_{728}}*{T_{27}}) +\phantom{111111111111111111111111}$\\
 &
& $+ 4374*({O_{364}}*{T_{27}}) + 19158152*({O_{546}}*{T_{27}})+\phantom{1111}$\\
 &
& $ + 12*({O_{182}}*{T_{27}})+ 818*({O_{24}}*{T_{27}}) +\phantom{111111111}$\\
 &
& $+3*({O_8}*{T_{27}})+ 8*({O_3}*{T_{27}}) + 3*({O_1}*{T_{27}}) $\\
\hline
85 & 2781283959047755505 & $ 278128395904775550578594372162961(O_{43046720}*T_3)+\phantom{11111111111}$\\
 & 78601345748937 & $+6973568964(O_{21523360}*T_3)+1377(O_{2532160}*T_3)+\phantom{}$\\
 & & $+15633(O_{223040}*T_3)+(O_{80}*T_3)+(O_1*T_3)$ \\
\hline
86 & 3992655605372987736 & $ 39926556053729877349446350778(O_{899590375372}*T_3)+\phantom{11111111111}$\\
 & 8420944965 & $+18974594184(O_{17299814911}*T_3)+3(O_1*T_3)$ \\
87 & 1438595418586772441 & $ 14385954185867724412000675312710(O_{832236432}*T_{27})+\phantom{11111111}$\\
 & 2000675395176& $+82465(O_{277412144}*T_{27})+(O_1*T_{27})$ \\
\hline
88 & 5474493800453436768 & $ 5474493800453436768255414641417732700(O_{59048}*T_3)+\phantom{1111111}$\\
 &266532979754636964 & $+11118288498762294(O_{29524}*T_3)+41841412812(O_{14762}*T_3)+\phantom{}$\\
 & & $+7780428000(O_{968}*T_3)+172898400(O_{484}*T_3)+\phantom{111111111}$\\
 & & $+43046721(O_{242}*T_3)+355758(O_{121}*T_3)+270(O_8*T_3)+\phantom{}$\\
 & & $+6(O_4*T_3)+3(O_1*T_3)$ \\
 89 & 5532420798784332770 & $ 5532420798784332769(O_{175289220588682799452640}*T_3)+(O_1*T_3)$ \\
\hline \hline
\end{tabular}

\newpage

\begin{tabular}{|r||r|r|}
\multicolumn{3}{c}{Table~2. Graph of operator $\Delta:
\mathbb{F}_3^n\to\mathbb{F}_3^n$}\\
\hline n &number of component  & connected components of  $G_n$ ($q=3$)\\
\hline\hline
90 & 6158701225597746845 & $ 615870122559774676380467939503795914(O_{720}*T_{19683})+\phantom{1111111}$\\
 & 86892156140567944 & $+8206424184863387028(O_{360}*T_{19683})+\phantom{1111111111111111111111}$\\
 & & $+31773262986(O_{240}*T_{19683})+119572(O_{120}*T_{19683})+\phantom{11111111}$\\
 & & $+243(O_{80}*T_{19683})+6(O_{40}*T_{19683})+2184(O_9*T_{19683})+\phantom{111}$\\
 & & $+8(O_3*T_{19683})+3(O_1*T_{19683})$ \\
\hline
91 & 119889609451754292363 & $ 11988960945175429076568214940160257228340(O_{728}*T_3)+\phantom{11111111}$\\
 & 20239673787574494747& $+159752024733625761705642(O_{364}*T_3)+1551804804(O_{182}*T_3)+\phantom{}$\\
 & & $+3720087(O_{104}*T_3)+5824(O_{91}*T_3)+9477(O_{56}*T_3)+\phantom{111111}$\\
 & & $+20412(O_{26}*T_3)+56(O_{13}*T_3)+104(O_7*T_3)+(O_1*T_3)$ \\
\hline
92 & 927094631514430173368 & $ 741675705206817198338502072873840(O_{31381059608}*T_3)+\phantom{11111111}$\\
 & 056790745879 & $+185418926295795674171523850821174(O_{15690529804}*T_3)+\phantom{1111111}$\\
 & & $+11817250826203332668808(O_{7845264902}*T_3)+\phantom{11111111111111111}$\\
 & & $+4251528(O_{3922632451}*T_3)+33354551686753680(O_{708584}*T_3)+\phantom{}$\\
 & & $+16677275843376837(O_{177146}*T_3)+6(O_{88573}*T_3)+3(O_8*T_3)+\phantom{}$\\
 & & $+3(O_1*T_3)$ \\
\hline
93 & 327025177475547816477 & $ 3270251774755478164778494348812300(O_{2668896516}*T_{27})+\phantom{11111}$\\
 & 8494349043735 & $+231434(O_{889632172}*T_{27})+(O_1*T_{27})$ \\
94 & 250315550501983041567 & $ 2503155505019830415674811918910117(O_{94143178826}*T_3)+\phantom{1111111}$\\
 & 4811918910126& $+6(O_{47071589413}*T_3)+3(O_1*T_3)$ \\
\hline
95 & 488457790885574358772 & $ 4651978960243204806224852640(O_{150094635296999120}*T_3)+\phantom{111111}$\\
 &4094360 & $+232598947411781699123246144(O_{37523658824249780}*T_3)+\phantom{11111}$\\
 & & $+1200757082375992968(O_{18761829412124890}*T_3)+\phantom{1111111111111}$\\
 & & $+16(O_{9380914706062445}*T_3)+2072(O_{14958320}*T_3)+\phantom{11111111}$\\
 & & $+518(O_{747916}*T_3)+(O_{80}*T_3)+(O_1*T_3)$ \\
\hline
96 & 119743402611383445461 & $ 11974340261053549305423610335949298861754(O_{19680}*T_{27})+\phantom{111}$\\
 & 19672660226266670629& $+84795240683752664742399630708(O_{9840}*T_{27})+\phantom{1111111111111}$\\
 & & $+94143178827(O_{6560}*T_{27})+28697814(O_{3280}*T_{27})+\phantom{111111111}$\\
 & & $+12309613114542720714(O_{240}*T_{27})+\phantom{1111111111111111111111}$\\
 & & $+46325417550228(O_{120}*T_{27})+177147(O_{80}*T_{27})+\phantom{111111111}$\\
 & & $+4374(O_{40}*T_{27})+435847140(O_{24}*T_{27})+1636(O_{12}*T_{27})+$\\
 & & $+270(O_8*T_{27})+6(O_4*T_{27})+8(O_3*T_{27})+3(O_1*T_{27})$\\
\hline
97 & 11612577085061309719 & $ 116125770850613097190293919508433(O_{54791330077120}*T_3)+\phantom{11111}$\\
 & 0312845199023& $+18925690589(O_{4214717698240}*T_3)+(O_1*T_3)$ \\
98 & 18620411870556350391 & $ 18620411870556350391681902671801614(O_{1025114613796}*T_3)+\phantom{1111}$\\
 & 682369545121455& $+466873315452(O_{512557306898}*T_3)+4374(O_{364}*T_3)+\phantom{11}$\\
 & & $+12(O_{182}*T_3)+3(O_1*T_3)$ \\
\hline
99 & 400732946193191571436 & $ 4007329461931915714365676385015741684388(O_{2178}*T_{19683})+\phantom{11}$\\
 & 8389252798466697422& $+2712867499128392424(O_{1089}*T_{19683})+\phantom{11111111111111111111}$\\
 & & $+283596580836(O_{726}*T_{19683})+39528(O_{363}*T_{19683})+\phantom{11}$\\
 & & $+243(O_{242}*T_{19683})+2(O_{121}*T_{19683})+(O_1*T_{19683})$ \\
\hline
100 & 49269609818912365316 & $ 49269609818912360530433137906892576934(O_{3486784400}*T_3)+\phantom{1111}$\\
 &500774994951451936 & $+4786067637088039181007(O_{2952400}*T_3)+\phantom{11111111111111111111}$\\
 & & $+5165694(O_{1476200}*T_3)+14528241(O_{80}*T_3)+54(O_{40}*T_3)+\phantom{}$\\
 & & $+3(O_8*T_3)+3(O_1*T_3)$ \\
\hline \hline
\end{tabular}
\end{center}

\section*{Appendix. Program on Mathematica with example.}

\begin{verbatim}
$HistoryLength=0;
<<Algebra`PolynomialPowerMod`

q = 3;(* The characteristic of field*)
res = {}; (*This is list of
lists :{length of sequence (this is n), indexTree,heightTree,
expansion of x^n -1 without (x-1)}*)

Do[poly = x^n -1;
  heightTree = q^(Length[NestWhileList[#/q &, n, Mod[#, q] == 0 &]] - 1);
  AppendTo[res, {n, indexTree = q^heightTree, heightTree,
  Factor[poly/(x - 1)^heightTree, Modulus -> q]}];
  , {n, 11, 13}];(*The region of sequence length*)

res // TableForm (* {n, indexTree,heightTree, expansion of x^i - 1
without (x - 1)} *)
\end{verbatim}

\(
\begin{array}{llll}
 11 & 3 & 1 & \left(2+2 x+x^2+2 x^3+x^5\right) \left(2+x^2+2 x^3+x^4+x^5\right) \\
 12 & 27 & 3 & (1+x)^3 \left(1+x^2\right)^3 \\
 13 & 3 & 1 & \left(2+2 x+x^3\right) \left(2+x^2+x^3\right) \left(2+x+x^2+x^3\right) \left(2+2 x+2 x^2+x^3\right)
\end{array}
\)

\begin{verbatim}
polynomials =    (*This is expansion of x^n - 1 : {polynom, its
degree, power of polynom}*)
    Table[
      If[res[[n, 3]] === res[[n, 1]], {{1, 0, 0}},
        If[Head[Last[res[[n]]]] ===Plus,
        {{temp = Last[res[[n]]], Length[CoefficientList[temp, x]] - 1, 1}},
          If[Head[Last[res[[n]]]] === Power, {{temp = First[Last[res[[n]]]],
                Length[CoefficientList[temp, x]] - 1, Last[Last[res[[n]]]]}},
            If[Head[Last[res[[n]]]] === Times &&  res[[n, 3]] > 1,
            {temp = First[#], Length[CoefficientList[temp, x]] - 1, Last[#]} &
            /@ (List @@ Last[res[[n]]]), {#, Length[CoefficientList[#, x]] - 1, 1} & /@
            (List @@ Last[res[[n]]])]]]],
 {n, Length[res]}];
 polynomials
\end{verbatim}

\noindent\(\left\{\left\{\left\{2+2 x+x^2+2 x^3+x^5,5,1\right\},
\left\{2+x^2+2 x^3+x^4+x^5,5,1\right\}\right\},\left\{\{1+x,1,3\},\left\{1+x^2,2,3\right\}\right\},\right.\\
\left.\left\{\left\{2+2x+x^3,3,1\right\},\left\{2+x^2+x^3,3,1\right\},\left\{2+x+x^2+x^3,3,1\right\},\left\{2+2
x+2 x^2+x^3,3,1\right\}\right\}\right\}\)

\begin{verbatim}
ord[polynom_, stepen_, qq_, jj_] :=(*The main procedure:
        the result is the order of (x - 1) in group by mod polynom^j*)
        (numbero = qq^(stepen(jj - 1))(qq^stepen - 1);
        (*stepen is degree of the polynom, qq == q*)
    polynomj = polynom^jj;
    dd = FactorInteger[numbero];

If[Times @@ (dd[[All, 2]] + 1) < 4000,
                     (*If[Length[divis]<4000,*)
  divis = Divisors[numbero];
  (*See article "Multiplicative Function Instead of Logarithm", Notation 2*)
  ii = 1;
  xm1degree[1] = PolynomialPowerMod[(x - 1), 1, {polynomj, qq}];
  While[Not[xm1degree[ii] === 1],
    kk = ii;
    ii++;
    While[Not[IntegerQ[divis[[ii]]/divis[[kk]]]], kk--];
    xm1degree[ii] =
      PolynomialPowerMod[xm1degree[kk],divis[[ii]]/divis[[kk]], {polynomj, qq}];];
  numbero/divis[[ii]] o[divis[[ii]]],

  order = numbero;
  (*See example in the Notation 2 of the article (about order of (x-1) by mod (x^2 +1)^2)*)
      Do[ll = dd[[ii, 2]];
          ff = dd[[ii, 1]];
          While[(temp = PolynomialPowerMod[x - 1, order, {polynomj, qq}]) === 1 && ll > 0,
              order /= ff; ll--];
          If[Not[temp === 1], order *= ff];
      ,{ii, Length[dd]}];
      numbero/order*o[order]]
);
\end{verbatim}

\begin{verbatim}
pluso[str_List, qq_] :=(*This is one multiplyer*)
  o[1] + Plus @@ Table[ord[str[[1]], str[[2]], qq, j], {j, str[[3]]}]
\end{verbatim}

\begin{verbatim}
timeso[arg_?MatrixQ, qq_] :=(*This is its product*)
  Expand[Times @@ (pluso[#, qq] & /@ arg)] //.
    {o[a_]^b_ -> a^(b - 1)o[a], o[a_]o[b_] -> GCD[a, b] o[LCM[a, b]]}
\end{verbatim}

\begin{verbatim}
r = timeso[#, q] & /@ polynomials                (*results*)
\end{verbatim}

\noindent\{o[1] + 2 o[121] + 243 o[242], 3 o[1] + 8 o[3] + 3 o[8]
+ 818 o[24],o[1] + 56 o[13] + 20412 o[26]\}

\begin{verbatim}
Table[{res[[i, 1]], Plus[r[[i]] /. o[a_] -> 1],
       Expand[T[res[[i,2]]]r[[i]]] /. {T[a_] -> Subscript[T, a],o[a_] -> Subscript[O, a]}},
{i,Length[r]}] // TableForm  (* results in TableForm *)
 \end{verbatim}

\noindent \(
\begin{array}{lll}
 11 \phantom{m}& 246\phantom{mmm} & O_1 T_3+2\ O_{121} T_3+243\ O_{242} T_3 \\
 12 & 832 & 3\ O_1 T_{27}+8\ O_3 T_{27}+3\ O_8 T_{27}+818\ O_{24} T_{27} \\
 13 & 20469 & O_1 T_3+56\ O_{13} T_3+20412\ O_{26} T_3
\end{array}
\)

\end{document}